\documentclass[12pt,centertags,reqno]{amsart}

\usepackage[foot]{amsaddr}
\usepackage{latexsym}
\usepackage[english]{babel}
\usepackage[T1]{fontenc}
\usepackage[numbers]{natbib}
\usepackage{amssymb}
\usepackage{fancyhdr}
\usepackage{url}
\usepackage{hyperref}
\usepackage{verbatim}
\usepackage{leftidx}
\usepackage{color,graphicx}
\usepackage{pdfpages}

\usepackage{amsmath}
\usepackage{algorithm}
\usepackage[noend]{algpseudocode}

\makeatletter
\def\BState{\State\hskip-\ALG@thistlm}
\makeatother

\usepackage{mathrsfs, mathtools}
\usepackage{stmaryrd}
\usepackage{nicefrac}
\usepackage{marvosym}
\mathtoolsset{showonlyrefs}

\numberwithin{equation}{section} \makeatletter
\renewcommand{\subsection}{\@startsection
	{subsection}{2}{0mm}{\baselineskip}{-0.25cm}
	{\normalfont\normalsize\bf}} \makeatother

\newtheorem{theorem}{Theorem}[section]

\theoremstyle{remark}
\newtheorem{remark}[theorem]{Remark}

\def \E {\mathcal E}
\def \F {\mathcal F}

\def \L {\mathcal L}

\def \P {\mathbf P}

\def \I {{\mathbf 1}}

\def \R {\mathbb R}
\def \bF {\mathbb F}

\def \bE {\mathbb E}

\newcommand{\ud}{\mathrm d}

\makeatletter
\def\thickhline{%
	\noalign{\ifnum0=`}\fi\hrule \@height \thickarrayrulewidth \futurelet
	\reserved@a\@xthickhline}
\def\@xthickhline{\ifx\reserved@a\thickhline
	\vskip\doublerulesep
	\vskip-\thickarrayrulewidth
	\fi
	\ifnum0=`{\fi}}
\makeatother

\newlength{\thickarrayrulewidth}
\begin{document}

	\author[R.~Frey]{R\"udiger Frey}\address{R\"udiger Frey, Institute for Statistics and Mathematics, Vienna University of Economics and Business, Welthandelsplatz, 1, 1020 Vienna, Austria }\email{rfrey@wu.ac.at}
	
	\author[V.~K\"ock]{Verena K\"ock}\address{Verena K\"ock, Institute for Statistics and Mathematics, Vienna University of Economics and Business, Welthandelsplatz, 1, 1020 Vienna, Austria }\email{verena.koeck@wu.ac.at}\thanks{We are grateful to Michaela Szoelgyenyi for useful remarks and suggestions}

	\title[Deep Neural Network Algorithms for Parabolic PIDEs]{Deep Neural Network Algorithms for Parabolic PIDEs and Applications in Insurance Mathematics}

	\date{}
	
	\begin{abstract}
	
  In recent years  a  large  literature on   deep learning based methods for the numerical solution partial differential equations has emerged; results for integro-differential equations on the other hand are scarce.   In this paper we study  deep neural network  algorithms for solving linear and semilinear parabolic partial integro-differential equations with boundary conditions in high dimension. To show the viability of our approach we discuss several case studies from insurance and finance.
		
	\end{abstract}

	\maketitle
	

\section{Introduction}

Many problems in  insurance and finance lead to terminal or boundary value problems involving parabolic partial integro-differential equations (PIDEs). Examples include option pricing in models with jumps,  the valuation of insurance contracts, ruin probabilities in non-life insurance,  optimal reinsurance problems and many applications in credit risk. These PIDEs can be linear (such as PIDEs arising in risk-neutral pricing) or semilinear (such as the dynamic programming equation in many stochastic control problems). Practical applications often involve several underlying assets or economic factors, so that one has to deal with  PIDEs in a high-dimensional space.     These PIDEs  do typically not admit an analytic solution, making the design of suitable  numerical methods an ongoing challenge.
Existing numerical  methods include deterministic schemes such as finite difference and finite element methods and random schemes based on Monte-Carlo methods. However, finite difference and finite element methods (see e.g. \citet{cont2005finite},  \citet{andersen2000jump}, \citet{matache2004fast}, \citet{kwon2011second}, \citet{briani2007implicit}) cannot be used in the case of high-dimensional PIDEs as they suffer from the {curse of dimensionality}.
Monte-Carlo  methods (see e.g. \citet{giles2008multilevel}, \citet{casella2011exact},  \citet{metwally2002using}) on the other hand are suitable for problems in higher dimensions. However, these methods  only provide a solution for a single fixed time-space point $(t,x)$.  This is problematic in risk management applications, where one needs to find the solution of a pricing problem for a large set $D$ of future scenarios.  The naive solution via nested Monte Carlo is  in most cases computationally infeasible.  Regression based  Monte Carlo methods (see \citet{bib:glasserman-03}) can sometimes help, but the choice of proper basis functions remains a delicate issue. Moreover, it is not straightforward to apply  Monte Carlo techniques  to  semilinear parabolic equations.

For these reasons many recent  contributions  study  machine learning techniques for the numerical solution of PDEs.    A large strand of this  literature is based on the representation of semilinear parabolic PDEs via backward stochastic differential equations (BSDEs). In the seminal papers \citet{han2017overcoming} and  \citet{bib:e-han-jentzen-17}, a BSDE is discretized using a time grid  $t_0 < t_1 < \dots < t_N=T$, and the solution at the initial date $t_0$ and its gradient  at  every time step $t_n$ are  approximated by a combined deep neural network. The parameters are trained by minimizing the difference between the network approximation and the  known terminal value of the solution.  An error estimation of this method is given by \citet{han2020convergence}.  \citet{kremsner2020deep} consider  an extension to  elliptic semilinear PDEs and applications to insurance mathematics.
Other  contributions  use a backward induction over the time steps $t_n$.  To begin with, \citet{hure2020deep}
estimate the solution and its gradient simultaneously by backward induction through sequential minimizations of suitable loss functions; moreover, they  provide convergence results for their method. The paper of \citet{beck2019deep} uses a different discretization method, called \emph{deep splitting}, that computes the unknown gradient of the solution by automatic differentiation, which reduces the size of the  networks. For linear PDEs one may use instead the simpler global regression approach of \citet{beck2018solving}. This paper   uses   a Feynman-Kac representation and the $\L^2$-minimality of conditional expectations to characterize  the solution  by means of an  infinite-dimensional  stochastic minimization problem which is solved with machine learning.
\citet{pham2021neural} combine the ideas of  \cite{hure2020deep}  and \cite{beck2019deep} to introduce a neural network scheme for fully nonlinear PDEs. Finally, \citet{germain2020deep} extend the method in \cite{hure2020deep} and  they provide a convergence analysis that can be adapted to show convergence of the deep splitting method of \cite{beck2019deep}.

Applications of deep learning methods to partial \emph{integro} differential equations on the other hand are scarce.   \citet{castro2021deep} presents an extension of \cite{hure2020deep} to PIDEs and he generalizes  the convergence results of \cite{hure2020deep} to his algorithm. Numerical case studies are however not provided. In fact, from a numerical viewpoint the method of \citet{castro2021deep} is   quite  involved, since one  needs to approximate the solution, the gradient and the  non-local term in the PIDE via three separate networks.     The work of  \citet{al2019applications}  is based on the \emph{deep Galerkin method} of \citet{sirignano2018dgm}. This is an alternative machine learning approach for PDEs, where the network is trained by directly minimizing the deviations of the fitted functions  from the desired differential operator and boundary conditions.

In this paper we consider  DNN algorithms for linear and semilinear parabolic PIDEs that generalize  the regression approach of \citet{beck2018solving}  and the deep splitting method of  \citet{beck2019deep}, respectively. In the semilinear case
we first linearize the equation locally in time using a time grid $t_n$, $n=0,\dots,N$.  Then we perform a backward induction over the grid points,  using in each step  the  DNN algorithm  for the linear case.  The advantage of this  approach, as opposed to the method of \citet{castro2021deep}, is the fact that we approximate only the solution by a DNN so that it suffices  to train a single network per time step.
In the semilinear case we propose a  alternative  linearization procedure to \citet{beck2019deep},  and we show in numerical experiments that this procedure performs substantially better than the original deep splitting  method of \citet{beck2019deep}.  Moreover, we  apply our DNN algorithms also to boundary value problems, while \citet{beck2018solving} and  \citet{beck2019deep} consider only   pure Cauchy problems.

The focus of our paper  is on applications to insurance and finance. Moreover, existing results on the convergence of DNN algorithms for PDEs give little guidance on how to construct an optimal network that achieves a given level of accuracy. For these reasons   an extension of  the convergence results from  \cite{hure2020deep} to PIDEs is left for future research. To assess accuracy and performance of our methodology  we  instead carry out extensive tests for  several multi-dimensional  PIDEs arising in actuarial mathematics. As a first test for the linear case we consider the
pricing of a stop-loss type reinsurance contract in a model where  claims arrive with stochastic intensity (see, e.g. \citet{grandell2012aspects}, \citet{ceci2020value}). In a second example we compute the ruin probability for an insurance company with three different business lines, which leads to a  boundary value problem. In both cases we  assess the performance of the deep learning method by comparing the solution to the results of an extensive   Monte Carlo simulation.  We go on and study  semilinear PIDEs arising in  stochastic control problems for jump diffusions.  The first test case is the multi-dimensional linear quadratic regulator problem, see for instance~\citet{oksendal2007applied}. The second case study is an optimization problem in  insurance. We consider an insurer who dynamically optimizes  his holdings of  some insurance portfolio in the presence of transaction costs and risk capital constraints.   In the absence of  capital constraints the problem admits  an analytic solution and  is therefore a useful test case. With capital constraints the model leads to a semilinear boundary value problem for which there is no explicit solution, and we study this case numerically.
Our numerical experiments show that in all test cases  the performance of the proposed approximation algorithms is quite satisfying in terms of accuracy and speed. 

The paper is organized as follows. Section \ref{sec:framework} introduces the general setting; Sections~\ref{sec:linearPIDE} and \ref{sec:examples-linear} deal with linear PIDEs,  and Sections~\ref{sec:semilinearPIDE} and~\ref{sec:examples-semilinear} are devoted to theory and case studies  concerning  the deep splitting method for semilinear PIDEs.

\section{Modeling Framework} \label{sec:framework}

We fix a probability space $(\Omega,\F, \P)$, a time horizon $T$ and a right continuous  filtration $\bF$. Consider measurable functions $\mu\colon [0,T] \times \R^d \to \R^d $, $\sigma\colon  [0,T] \times \R^d \to \R^{d\times d}$ and $\gamma^X \colon [0,T] \times \R^d \times E \to \R^d$, where $(E,\E)$ is a separable measurable space.
We assume that $(\Omega,\F, \P)$  supports  a $d$-dimensional Brownian motion $W$   and a Poisson random measure $J$ on $[0,T]\times E$. The compensator of $J$ is given by  $\nu(\ud z) \ud t$ for a sigma-finite measure $\nu$ on $E$.  We consider a $d$-dimensional process  $X$ that is  the unique strong solution to the SDE
\begin{align} \label{eq:processX}
	\ud X_t= \mu(t,X_t)  \ud t + \sigma(t,X_t)  \ud W_t +  \int_{E} \gamma^X(t,X_{t-},z)  J(\ud t, \ud z), \qquad X_0=x \in \R^d.
\end{align}
Define the set
$D^+(t,x) = \{ z \in E \colon \|\gamma^X(t,x,z)\| >0 \}$ and note that $X$ jumps at $t$ whenever $J\big (\{t\} \times D^+(t,X_{t-}) \big ) >0$. We assume that
\begin{equation}\label{eq:finite-activity}
\bE \left [ \int_0^T \nu(D^+(t,X_{t-})) \ud t \right ] < \infty\,.
\end{equation}
This condition ensures  that $X$ has a.s.~only finitely many jumps on $[0,T]$, so that every integral with respect to $J$ is well-defined. The restriction  to finite activity jump processes simplifies the exposition and it is sufficient for most applications in insurance.  Conditions on the coefficients $\mu$, $\sigma$ and $\gamma^X$   ensuring that the SDE \eqref{eq:processX} has a unique strong solution are given for instance in \citet{bib:gihman-skohorod-80} or in \citet{bib:kliemann-koch-marchetti-90}.

By $\mathcal{C}^k(D) $ we denote the functions that are $k$ times continuously differentiable on the set $D \subset\R^d$, and by $\mathcal{C}^{1,k}([0,T] \times D)$  we denote functions that
are once continuously differentiable in $t$ and $k$-times continuously differentiable  in $x$ on the set $[0,T] \times D$. For every function $h\in \mathcal{C}^{1,k}([0,T] \times D)$  we write $h_{x_i}$ for the first derivatives of $h$ with respect to $x_i$ for $i \in \{1,\dots, d\}$ respectively, $h_{x_i x_j}$ for second derivatives, for $i,j \in \{1,\dots, d\}$, and finally $h_t$ denotes the first derivative with respect to time.

Define the matrix $\Sigma(t,x)=(b_{i,j}(t,x), \ i,j=1,\dots,d)$, by
$$\Sigma(t,x)=\sigma(t,x)\sigma^\top(t,x)   $$
and consider for $u\in \mathcal{C}^{1,2}([0,T]\times \R^d)$ the  integro-differential operator $\mathcal{L}$ given by
\begin{align}\label{eq:generatorX}
\L u(t,x) :=&\sum_{i=1}^{d} \mu_i(t,x) { u_{x_i}}(t,x) +\frac{1}{2} \sum_{i,j=1}^{d} b _{i,j}(t,x)  u_{x_i x_j}(t,x) \\ &+ \int_{\R^d} [u(t,x+\gamma^X(t,x,z))-u(t,x)] \nu(\ud z), \quad x \in \R^d, \ t \in [0,T]\,.
\end{align}
The operator $\mathcal{L}$ is the generator of $X$, and $X$ is a solution of the martingale problem for $\L$, see \citet{bib:ethier-kurtz-86} for details.

Consider  functions  $c\colon [0,T] \times \R^d \to \R$, $r\colon [0,T]\times \R^d \to \R$ and  $g \colon [0,T] \times \R^d \to \R$, and let  $D$ be an open subset of $\R^d$.  In Section~\ref{sec:linearPIDE} we are interested in the following boundary value  problem
\begin{align} \label{eq:Cauchyproblem}
&u_t(t,x)+\L u (t,x)-  r(t,x)u(t,x) + c(t,x)=0, \quad (t,x) \in [0,T) \times D,\\
&u(t,x)=g(t,x), \quad (t,x) \in \big([0,T) \times (\R^d \setminus D)\big)  \cup \big ( \{T\}\times \R^d \big ).
\end{align}
The special case $D = \R^d$ corresponds to a pure Cauchy problem without boundary conditions; in that case we use the simpler notation $g(T,x) =: \varphi(x) $ to denote  the terminal condition.

It follows from the Feynman-Kac formula that  under some integrability conditions a classical solution $u$ of \eqref{eq:Cauchyproblem} has the probabilistic representation
\begin{equation} \label{eq:linearsolution}
u(t,x)=	\bE \Big[\int_{t}^{T\wedge \tau}\hspace{-0.2cm} e^{- \int_{t}^{s}r(u, X_u)  \ud u}c(s, X_s) \ud s+e^{- \int_{t}^{T\wedge \tau}r(u, X_u)  \ud u}g \big(T\wedge \tau, X_{T\wedge \tau}\big) \mid  X_t=x\Big], \;
\end{equation}
where $\tau:=\inf \{ s\ge t \colon  X_s \notin D \}$. Conversely, \citet{bib:pham-98} and \citet{bib:colaneri-frey-21} provide conditions ensuring that the function  $u$ defined in \eqref{eq:linearsolution} is a classical solution of the problem \eqref{eq:Cauchyproblem} (for the case of a pure Cauchy problem). Further existence results for linear PIDEs include \citet{bib:gihman-skohorod-80}, \citet{bensoussan1982optimal}, and \citet{davis2013jump}.

In Section~\ref{sec:linearPIDE} we propose a deep neural network (DNN) algorithm  to approximate the function $u$ defined in \eqref{eq:linearsolution}. In Section \ref{sec:semilinearPIDE} we are interested in semilinear problems of the form
\begin{align} \label{eq:semilinear_problem_2}
&u_t(t,x)+\L u (t,x) + f(t,x,u(t,x),\nabla u(t,x)) =0 , \quad (t,x) \in [0,T) \times D,\\
&u(t,x)=g(t,x), \quad (t,x) \in \big([0,T) \times (\R^d \setminus D)\big) \cup  \big(\{T\}\times \R^d\big),
\end{align}
where $f \colon [0,T]\times \R^d \times \R \times \R^d \to \R$ is a nonlinear function such as the Hamiltonian in a typical Hamilton Jacobi Bellman equation. To handle this case we partition  the interval $[0,T]$ into time points $0=t_0 < t_1 <\dots <t_N =T$, consider a linearized version of \eqref{eq:semilinear_problem_2} for each subinterval $[t_n, t_{n+1}]$,  and apply the DNN algorithm  for the linear case recursively.

\section{Deep neural network  approximation for linear PIDEs} \label{sec:linearPIDE}

In this section we consider linear PIDEs. We extend the regression-based algorithm of \citet{beck2018solving} to PIDEs  and we include boundary conditions into the analysis.

\subsection{Representation as solution of a minimization problem}

Fix some time point $t \in [0,T)$ and a closed and bounded set $A \subset \overline D$. Define  the  function $u\colon  [0,T] \times \R^d \to \R$  by the Feynman-Kac representation~\eqref{eq:linearsolution}. We want to compute an approximation to the function  $ u(t,\cdot)$ on the set $A$.  The key idea is to write this function   as solution of a minimization problem on an infinite dimensional space. This representation is used to construct a loss function for our deep neural network  method.

Consider some random variable $\xi$ whose distribution is absolutely continuous with respect to the Lebesgue measure such that the corresponding density has support $A$ (in applications the distribution of $\xi$ is often the uniform distribution on $A$) and denote by $X^{\xi}$ the solution of the SDE \eqref{eq:processX} with initial value $X_t = \xi$.
Define the random variable
		\begin{align*}
			Y^\xi:=\int_{t}^{T\wedge \tau} e^{- \int_{t}^{s}r(u, X_u^\xi)  \ud u}c(s, X_s^\xi) \ud s+e^{- \int_{t}^{T\wedge \tau}r(u, X_u^\xi)  \ud u}g \Big(T\wedge \tau, X_{T\wedge \tau}^\xi\Big)
		\end{align*}
Assume that  $\bE[|Y^\xi|^2]<\infty$ and that the function $u(t,\cdot)$ belongs to $\mathcal{C}^0(\overline A)$.  Since  $X^\xi$ is a Markov process it  holds that  $u(t,\xi) = \bE [Y^\xi \mid \sigma(\xi)]$, where $\sigma(\xi)$ is the sigma-field generated by $\xi$.   Since $Y^\xi$ is square integrable we thus get from the $\mathcal{L}^2$-minimality of conditional expectations that
\begin{align}\label{eq:minprob}
			\bE &\Big[\big|Y^\xi- u(t,\xi) \big|^2 \Big] = \inf \left \{ \bE \Big[\big|Y^\xi- Z \big|^2 \Big] \colon Z \in L^2 (\Omega, \sigma(\xi),\P) \right \}\,.
\end{align}
Since $u(t,\cdot) \in \mathcal{C}^0({A})$  and since the density of $\xi$ is strictly positive on $A$  we conclude that $u(t,\cdot)$ is  the unique solution of the minimization problem
\begin{equation}\label{eq:minimization-for-u}
\min \bE \Big[\big|Y^\xi- v(\xi) \big|^2 \Big]\,, \quad  v \in \mathcal{C}^0({A})\,.
\end{equation}
The problem~\eqref{eq:minimization-for-u} can be solved with deep learning methods, as we explain next.

\subsection{The algorithm}
The first step in  solving \eqref{eq:minimization-for-u} with machine learning techniques is to simulate trajectories of ${X^\xi}$ up to the stopping time $\tau$. The simplest method is  the Euler-Maruyama scheme. Here we choose a time discretization $t=t_0<t_1<\dots<t_M=T$, $\Delta t_m = t_{m}-t_{m-1}$,\footnote{%
We use $m$ to index the time steps in the Euler-Maruyama scheme and $n$ to index the grid points used in the linearization step of the deep splitting method in Section~\ref{sec:semilinearPIDE}.}
generate $K$ simulations $\xi^{(1)},\dots,\xi^{(K)}$ of the random variable $\xi$ and simulate $K$ paths $X^{(1)},\dots,X^{(K)}$ of $X^\xi$ up to the stopping time $\tau$ by the following recursive algorithm. We let $ X_{t}^{(k)} =\xi^{(k)}$, and for $m \ge 1$,
	\begin{align}\label{eq:EuMa1}
		{X}_{t_{m} \wedge \tau}^{(k)} := 	{X}_{t_{m-1} \wedge \tau }^{(k)} + \I_{(0,\tau)} (t_{m-1}) \ \Big( &\mu({t_{m-1}},{X}^{(k)}_{t_{m-1}}) \Delta t_m +  \sigma(t_{m-1},{X}^{(k)}_{t_{m-1}}) \  \Delta W^{(k)}_{t_m} \\& \quad +
\int_{t_{m-1}}^{t_{m}} \int_{\R^d} \gamma(t_{m-1},X_{t_{m-1}}^{(k)},z) \ J^{(k)}(\ud z, \ud s) \Big)\,.
	\end{align}
Note that the integrand in the integral with respect to $J^{(k)}$ is evaluated at $t_{m-1}$ so that this integral corresponds to the increment of a standard compound Poisson process. Using these simulations we compute for each path
$$ Y^{(k)}:=\int_{t}^{T\wedge \tau} e^{- \int_{t}^{s}r(u, X^{(k)}_u)  \ud u}c(s, X^{(k)}_s) \ud  s+e^{- \int_{t}^{T\wedge \tau}r(u, X^{(k)}_u)  \ud u} g({T \wedge \tau}, X^{(k)}_{T\wedge \tau}),$$
where the integrals on the right can be approximated by Riemann sums.

	In the next step we approximate $u(t,\cdot )$ by a deep neural network $\mathcal{U}_t(\cdot) = \mathcal{U}_t(\cdot; \theta) \colon A \to \R^d$. We determine the network parameters $\theta$ (training of the network) by minimizing the  loss function
	\begin{align} \label{eq:exploss}
		\theta \mapsto \frac{1}{K}\sum_{k=1}^{K} \big ( Y^{(k)}-\mathcal{U}_t(\xi^{(k)};\theta)\big)^2\,.
	\end{align}
For this we rely on  stochastic gradient-descent methods; algorithmic details are given in the next section.  This approach  can be considered  as a \emph{regression-based} scheme since one attempts  to minimize the squared error between the DNN approximation $\mathcal{U}_t(\cdot;\theta)$ and the given terminal and boundary values of the PIDE.

\section{Examples for the linear case: } \label{sec:examples-linear}

In this section we test the performance of the proposed DNN algorithm for linear PIDEs in two  case studies.   First we price a reinsurance contract in  the  model of~\citet{ceci2020value},  where the claims process follows a doubly stochastic risk process. In the second example  we compute the ruin probability of a non-life  insurance company  with several business lines, which leads to a boundary value problem.

\subsection{Valuation of an insurance contract with doubly stochastic Poisson arrivals}	
\begin{figure}[h]
	\centering
	\begin{tabular}{@{}c@{\hspace{.5cm}}c@{}}
		\includegraphics[width=0.95\textwidth]{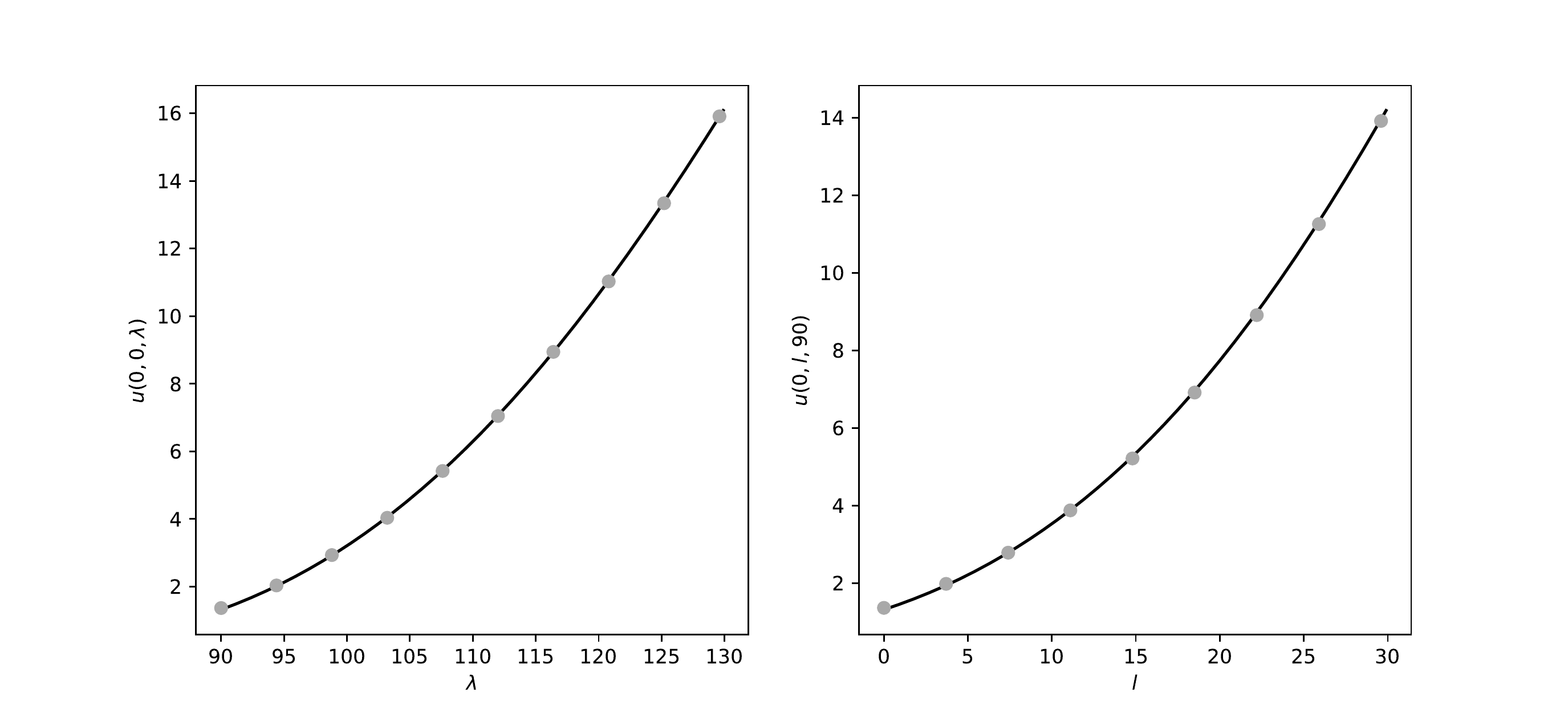}
	\end{tabular}
	\caption{ Solution $u(0,l,\lambda)$ for  $l=0$, $\lambda \in [90,130]$ (\textit{left}) and $l \in [0,30]$, $\lambda=90$ (\textit{right}) computed with the DNN-algorithm (\textit{black line}) and reference points computed with Monte-Carlo (\textit{grey dots}).}
	\label{fig:losslayer}
\end{figure}

We consider an insurance company and a reinsurer who enter into a reinsurance contract with a given maturity $T=1$.
To model the losses in the insurance portfolio underlying this contract we consider a sequence $\{T_n\}_{n \in \mathbb N}$ of claim arrival times with nonnegative intensity process  $\lambda^L=(\lambda^L_t)_{t \ge 0}$  	
and a sequence $\{Z_n\}_ {n\in \mathbb N}$ of claim sizes that are iid strictly positive random variables independent of the counting process $N=(N_t)_{t \geq 0}$ defined by $N_t = \sum_{n=1}^{\infty} \I_{\{T_n \le t\}}$. The loss process $L=(L_t)_{t \geq 0}$ is given by
$
L_t=\sum_{n=1}^{N_t}Z_n .
$
We assume that the $Z_n$ are  Gamma($\alpha$,$\beta$) distributed with density $f_{\alpha,\beta}(z)$.  This is a common choice in insurance. Moreover, the Gamma distribution is closed under convolution so that the sum of independent Gamma distributed random variables can be generated with a single simulation, which speeds up the sampling of trajectories from  $L$. The claim-arrival intensity process $\lambda^L$ satisfies the  SDE
\begin{align}\label{eq:process_lambda}
\ud \lambda^L_t&=  b( \lambda^L_t) \ud t + \sigma(\lambda^L_t) \ud W_t, \quad \lambda^L_0 =\lambda_0\in \R_+,
\end{align}		
where $W$ is a standard Brownian motion. In this example it is convenient to write the process $X$ in the form $X_t = (L_t,\lambda^L_t)$.
We assume that the reinsurance contract is a \textit{stop-loss} contract, i.e. the indemnity payment  is of the form $\varphi(L_T)$ with
\begin{align} \label{eq:maxplus}
\varphi(l) = [l-{K}]^+, \quad  \text{with }  [z]^+ = \max\{z,0\} .
\end{align}		
The {market value} $u$ at time $t \in [0,T]$ of the reinsurance contract is defined  by
\begin{align}\label{eq:market-value}
u(t,l,\lambda) :=\bE [ \varphi(L_T)|L_t=l,\lambda^L_t=\lambda], \quad (l,\lambda) \in \R^0_+ \times \R_+
\end{align}	
\citet{ceci2020value} show  that $u$ is the unique solution of the  PIDE
$
u_t(t,l,\lambda)+ \mathcal L u(t,l,\lambda)=0
$
with terminal  condition $u(T,l, \lambda)=\varphi(l)$ and generator
\begin{align*}
\mathcal L u(t,l,\lambda) = u_\lambda(t,l,\lambda) b(\lambda) + \frac{1}{2}u_{\lambda \lambda} (t,l,\lambda) \sigma(\lambda)^2 + \lambda \int_{\R} [u(t,l+z,\lambda)-u(t,l,\lambda)] f_{\alpha,\beta}(z) \ud z\, ,
\end{align*}
for $(l,\lambda) \in \R^0_+ \times \R_+$, $t\in[0,T)$. There is no explicit solution for this PIDE, and  we approximate $u(0,l,\lambda)$ on the set $A := \{ (l,\lambda) \colon l\in [0,30], \lambda \in [90,130]\}$ with a deep neural network $\mathcal{U}_0(l,\lambda)$. Table \ref{tab:parameters_doubly} contains the parameters we use. Paths of the processes $L$ and $\lambda^L$ are simulated with the Euler-Maruyama scheme and $\xi \sim \text{Unif}(A)$.
\begin{table}[h]
	\begin{center}
		\begin{tabular}{c c c c c c  }
			$b(\lambda)$& $\sigma(\lambda) $ &  $\alpha$ & $\beta$& $K$ & \\
			\hline
			$ 0.5  (100-\lambda)$& $0.2  \lambda $ &  $1$ & $1$& $90$ & \\
			\hline
		\end{tabular}
	\end{center}
	\vspace{2mm}
	\caption{Parameters used for the valuation of the stop-loss contract. The claim sizes are Gamma$(\alpha,\beta)$ distributed with density function $f_{\alpha,\beta}$.}  \label{tab:parameters_doubly}
\end{table}			
Figure \ref{fig:losslayer}
shows the approximate solution $\mathcal{U}_0$ obtained by the DNN algorithm; details on the training procedure are given in Remark~\ref{rem:num_detail_linear}.

As a reference we compute for fixed $(l, \lambda)$ approximate values  $U^{MC}(l,\lambda) \approx u(0,l,\lambda)$ with Monte-Carlo using $10^6$ simulated paths for each point $(l, \lambda)$  (paths are simulated with the Euler-Maruyama scheme).
The \emph{relative $L^1$-error} between the DNN approximation $\mathcal{U}$ and the MC-solution $U^{MC}$ is defined as
\begin{align}\label{eq:relative-error}
\epsilon:=\mathbb{E} \Big[
\Big|\frac{\mathcal{U}(\xi^l, \xi^\lambda)-{U}^{MC}(\xi^l, \xi^\lambda)}{{U}^{MC}(\xi^l, \xi^\lambda)}\Big|\Big].
\end{align}
Using 1000 simulations of $\xi^l \sim \text{Unif}([0,30]), \xi^\lambda \sim \text{Unif}([90,130])$ we obtained a relative error of {$\epsilon=0.0018$}. On a Lenovo Thinkpad notebook with an Intel Core i5 processor (1.7 GHz) and 16 GB of memory the computation of $\mathcal{U}$ via the training of a DNN  took around 322 seconds, whereas the computation of $U^{MC}(l,\lambda)$  with Monte-Carlo for a fixed point $(l,\lambda)$ took around 4.3 seconds. This shows that the DNN approach is faster than the MC approach if one wants to compute $u(t,l_i,\lambda_i)$ for a grid $(l_i,\lambda_i)$ with more than 100 grid points. Note also that the training of the networks can be accelerated further by using  GPU computing.

\begin{remark}[Details regarding the numerical implementation] \label{rem:num_detail_linear}		
Our choice of the network architecture largely follows  \cite{beck2018solving}. Throughout we use a  DNN with 2 hidden layers. In the reinsurance example we worked with  50 nodes for the hidden layers. For a payoff of the form  \eqref{eq:maxplus} it is advantageous to use  \texttt{softplus} as activation function for the hidden layers and the \texttt{identity} for the output layer. The networks are initialized  with random numbers using the Xavier initialization. We use mini-batch optimization with mini-batch size $M=5000$, batch normalization and 10000 epochs for training. We minimize the loss function with the Adam optimizer and decrease learning rate from 0.1 to 0.01, 0.001 and 0.0001 after a 2000, 4000 and 6000 epochs.
\end{remark}

\subsection{Ruin probability of an insurer with different business lines}

In the second case study we consider a non-life insurer who provides financial protection against losses from earthquakes, storms and floods over a given time horizon $T=1$. The occurrence of these natural disasters is modeled using independent Poisson counting processes $N^E, N^S,  N^F$  with constant intensities $\lambda^E, \lambda^S,  \lambda^F$. The magnitudes of hazards caused by earthquakes,  storms or  floods are modeled by iid sequences $\{E_n\}_ {n\in \mathbb N}$, $\{S_n\}_ {n\in \mathbb N}$, $\{F_n\}_ {n\in \mathbb N}$ of  positive random variables.

Furthermore we consider  catastrophic events that cause simultaneous losses in the  windstorm and flood business lines. These events  follow a Poisson process $N^{S,F}$ with constant intensity $\lambda^{S,F}$ and iid claim sizes sequences  $\{\widetilde S_n\}_ {n\in \mathbb N}$, $\{\widetilde F_n\}_ {n\in \mathbb N}$. Define $\lambda=(\lambda^E, \lambda^S,  \lambda^F, \lambda^{S,F}) \in \R^5_+$.  Then the risk process $X=(X_t^E,X_t^S,X_t^F)_{t\geq 0}$ is given by
\begin{align} \label{eq:loss_CAT}
X^E_t&=b^E t+ \sigma^E W_t^E-\sum_{n=1}^{N^E_t}E_n\,, \quad t \geq 0, \qquad \textit{(earthquake)} \\
X^S_t&=b^S t+\sigma^S W_t^S-\sum_{n=1}^{N^{S}_t} S_n-\sum_{n=1}^{N^{S,F}_t} \widetilde S_n\,, \quad t \geq 0, \qquad \textit{(storm)} \\
X^F_t&=b^F t+\sigma^F W_t^F-\sum_{n=1}^{N^{F}_t}F_n-\sum_{n=1}^{N^{S,F}_t} \widetilde F_n, \quad t \geq 0, \qquad \textit{(flood)}
\end{align}
where $W=(W_t^E,W_t^S,W_t^F)_{t\ge 0}$ is a $3$-dimensional standard Brownian motion, $\sigma=(\sigma^E,\sigma^S,\sigma^F) \in \R^3_+$  and $b=(b^E,b^S,b^F) \in \R^3$. The continuous terms $b^i t+\sigma^i W^i_t, \, i \in \{E,S,F\} $ are an  approximation of the difference of cumulative premium payments and small losses  up to time $t$ for each business line. We assume that all claim sizes $E, S, F, \widetilde S, \widetilde F$ are Gamma($\alpha_i$,$\beta_i$) distributed for $i \in I:= \{1,2,3,4,5 \}$. The generator of $X$ is  given by
\begin{align*}
\mathcal Lu(t,x) = \sum_{i=1}^3 u_{x_i}(t,x) b^i	+\frac{1}{2}\sum_{i=1}^3u_{x_ ix_i}(t,x)(\sigma^i)^2 + \int_{\R_+^3} [u(t,x-z)-u(t,x)] \nu(\ud z), 	\end{align*}	
where
$\nu(\ud z)=\sum_{i =1}^3 \lambda^i f_{\alpha_i,\beta_i}(z_i) \ud z_i + \lambda^4 f_{\alpha_{4},\beta_{4}}(z_2)f_{\alpha_{5},\beta_{5}}(z_3) \ud z_2 \ud z_3. $

We define the ruin time $\tau<T$ as the first time the minimum of all three risk processes falls below zero, i.e. $\tau:=\inf \{t \ge 0 \colon  \min(X^E_t,X^S_t,X^F_t)\leq 0\}$. Define the set
$D:= (0,\infty)^3 \subset \R^3$ and note that $\tau$ is the first exit time of the risk process $X$ from $D$. Put
\begin{align*}
g(t,x)=
\begin{cases}
0, \quad t=T, x \in D, \\
1, \quad t\leq T, x \notin D.
\end{cases}
\end{align*}
The ruin probability $u$ at time $t \in [0,T]$  given $\tau > t$ is
\begin{align}\label{eq:market-value2}
u(t,x) := \P(\tau \le T \mid X_t=x) =  \bE 	[ g(T \wedge \tau ,X_{T\wedge \tau})|X_t=x], \quad x \in D 	.
\end{align}	
By standard arguments $u$ is the unique solution of the linear PIDE
\begin{align}\label{eq:backward2}
&u_t(t,x)+ \mathcal Lu(t,x)=0,\quad (t,x)\in [0,T)\times D,
\end{align}
with boundary condition $v(t,x)=g(t,x)$.

\begin{figure}[h]
	\centering
	\begin{tabular}{@{}c@{\hspace{.5cm}}c@{}}
		\includegraphics[width=0.95\textwidth]{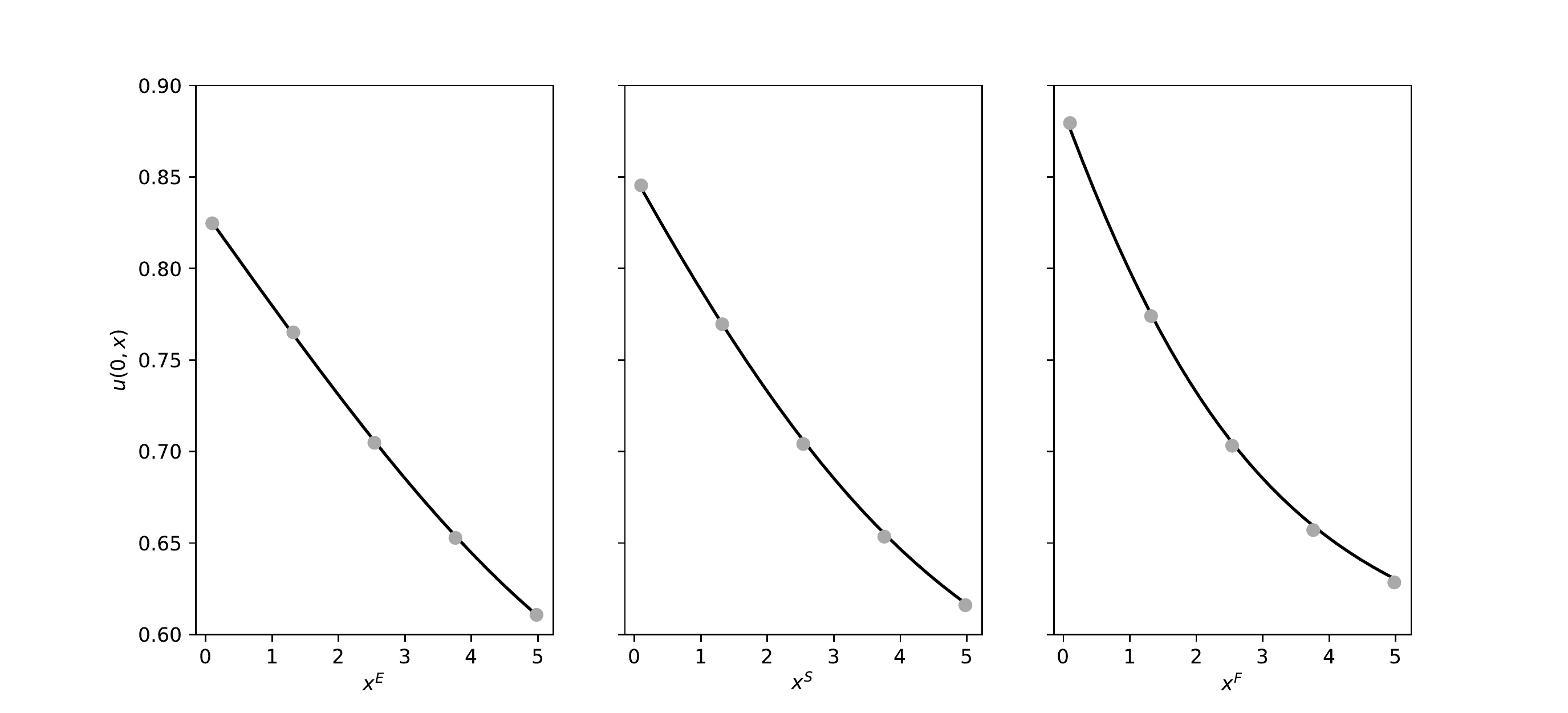}
	\end{tabular}
	\caption{ Solution $u(0,x)$ on  $\{(x^E,3,3) \colon  x^E \in [0.1,5] \}$ (\textit{left}), $\{(3,x^S,3) \colon  x^S \in [0.1,5] \} $ (\textit{middle}), and $\{(3,3,x^F) \colon  x^F \in [0.1,5] \} $ (\textit{right}) with the  DNN-algorithm (\textit{black line}) and Monte-Carlo (\textit{grey dots}).}
	\label{fig:cat}
\end{figure}

We approximate $x \mapsto u(0,x)$ on $A:=[0.1,5]^3$ with a deep neural network $\mathcal{U}_0$ assuming that $\xi \sim \text{Unif}(A)$. Moreover we  compare the result to Monte-Carlo approximations $U^{MC}(x)$ for fixed points $x \in A$.  Table \ref{tab:parameters_CAT} contains the parameters we use.
\begin{table}[h]
	\begin{center}
		\begin{tabular}{c c c c c    }
			$\sigma$   & 	$b$ & $\lambda $  &   $\alpha $ & $\beta$
			\\ \hline
			\rule{0pt}{0.9\normalbaselineskip}
			$(0.1,0.1,0.1)$  &  $(6,6,6)$  &$(2,2,10,1)$ &$(3,2,0.5,2,1)$ & $(1,1,1,1,1)$			
		\end{tabular}
	\end{center}
	\vspace{2mm}
	\caption{Parameters used for the case study on  ruin probabilities. Recall that claim sizes are Gamma$(\alpha_i,\beta_i)$ distributed.}  \label{tab:parameters_CAT}
\end{table}	
Figure \ref{fig:cat} shows the approximate solution $\mathcal{U}_0(x) \approx u(0,x)$ obtained by the DNN algorithm on  the sections $\{(x^E,3,3) \colon  x^E \in [0.1,5] \}$, $\{(3,x^S,3) \colon  x^S \in [0.1,5] \} $ and $\{(3,3,x^F) \colon  x^F \in [0.1,5] \} $. The network architecture is as in Remark \ref{rem:num_detail_linear}, but since we are working in  $d=3$ dimensions we choose mini-batch size $M=6000$ and 100 nodes for the hidden layers. 	To verify our result we computed $U^{MC}(x)\approx u(0,x)$ with Monte-Carlo for fixed  $x \in A$ using $2\cdot10^6$ path simulations for each point. The relative error \eqref{eq:relative-error} was computed to $\epsilon=0.0016$ (using 1000 simulations of $\xi \sim \text{Unif}(A)$), which is again very small.  Training of the network took  around 740 seconds; the computation of the reference solution via Monte Carlo took around 20 seconds per point.

This example clearly shows the advantages of the DNN method over standard Monte Carlo for computing the solution on the entire set  $A$. Suppose that we want to compute the solution on $[0.1, 5]^3$ (as in our case study). Even the very coarse grid $\{1,2,\dots,5\}^3$ has already 125 gridpoints, and  computing the solution for each gridpoint takes approximately $20 \times 125 = 2500$ seconds, which is about three times the time for training the network.

\section{Deep learning approximation for semilinear PIDEs}\label{sec:semilinearPIDE}

Next  we consider   semilinear PIDEs of the form
\begin{align} \label{eq:semilinear_problem}
u_t(t,x) &+\L u (t,x) + f(t,x,u(t,x),\nabla u(t,x)) =0 , \quad (t,x) \in [0,T) \times D,\\
 u(t,x)  &=g(t,x), \quad (t,x) \in \big([0,T) \times (\R^d \setminus D)\big) \cup  \big(\{T\}\times \R^d \big)\,.
\end{align}
Here $f\colon [0,T]\times \R^d \times \R \times \R^d \to \R$ is a nonlinear function and the integro-differential operator $\L$ is given by \eqref{eq:generatorX}.
Our goal is to  extend the deep splitting approximation method developed by \citet{beck2019deep} for semilinear PDEs without boundary conditions to the more general equation~\eqref{eq:semilinear_problem}. Moreover,  we propose an approach for  improving the performance of the method. In the derivation of the deep splitting method we implicitly assume that a classical solution of this PIDE exists; see  for instance~\citet{davis2013jump} for results on this issue.

\subsection{Basic method}
We  divide the interval $[0,T ]$ into subintervals using a  time grid $0=t_0<t_1<\dots<t_N=T$ and we let  $\Delta t_n = t_{n}-t_{n-1}$.
Suppose that $u^*$ is the unique  solution of the PIDE~\eqref{eq:semilinear_problem}. Define for a fixed $u\in \mathcal{C}^{0,1}([0,T] \times D)$  the function $$(t,x) \mapsto f_{[u]}(t,x)=f(t,x,u(t,x),\nabla u(t,x)).$$
	Then $u^*$ also solves   the linear PIDE
	\begin{align} \label{eq:tempdis}
		u_t(t,x)+\L u (t,x) + f_{[u^*]}(t,x) =0.
	\end{align}		
 Applying the Feynman-Kac formula  to the linear  PIDE~\eqref{eq:tempdis} we get the following representation for $u^*$
\begin{align} \label{eq:splitting_PIDE}
		{u}^* (t_{n-1},x)=\mathbb{E} \Big[ u^*(t_{n}\wedge \tau, X^x_{t_{n}\wedge \tau})+ \int_{t_{n-1}}^{t_{n}\wedge \tau} f_{[u^*]}(s,X^x_s) \ud s  \Big ], \quad x \in D. 	
\end{align}
	Here $X^x$ solves the SDE \eqref{eq:processX} with initial condition $X_{t_{n-1}}^x =x$, that is
	\begin{align*}
	X^x_t=x+\int_{t_{n-1}}^{t} \mu(s,X_s^x) \ \ud s + \int_{t_{n-1}}^{t}\sigma(s,X_s^x) \ \ud W_s +  \int_{t_{n-1}}^{t} \int_{\R^d} \gamma(s,X_{s-}^x,z) \ J(\ud z, \ud s),	
	\end{align*}
and the stopping time $\tau$ is given by $\tau = \inf \{s > t_{n-1} \colon X^x_s \notin D\}$.
	In the  deep splitting method proposed in \cite{beck2019deep}  the integral term  in \eqref{eq:splitting_PIDE} is approximated  by the \emph{endpoint rule}
	\begin{align}\label{eq:endpoint1}
	\int_{t_{n-1}}^{t_{n}} f_{[u^*]}(s,X^x_s) \ud s \approx  f_{[u^*]}(t_{n},X^x_{t_{n}}) \Delta t_n\,.
	\end{align}
 An approximate solution $\mathcal{U}_{n}(\cdot)  \approx u^*(t_n,\cdot)$, $n= 0,\dots N$ is then found by backward induction over the grid points. The solution at the maturity $t_N=T$ is defined using the terminal condition, that is
	$	{\mathcal U}_N (x):=g(T,x).	$
Given $\mathcal{U}_{n}$ we then compute $\mathcal{U}_{n-1}$, $n=1,2,\dots,N$  as follows. First we put
	\begin{align} \label{eq:basic-linearization}
	\tilde{u}(t_{n-1},x) &:=  \mathbb{E} \Big[{\mathcal U}_n(X^x_{t_n }) \I_{\{\tau >t_n\}} + g(\tau, X^x_\tau)   \I_{\{\tau \le t_n\}} \\&+  (t_n \wedge \tau - t_{n-1})  f\big (t_n \wedge \tau ,  X^x_{t_n \wedge \tau}, {\mathcal U}_n ( X^x_{t_n \wedge \tau}), \nabla {\mathcal U}_n ( X^x_{t_n \wedge \tau})  \big)  \Big]\,.
	\end{align}
	We then compute  $\mathcal{U}_{n-1} \approx \tilde{u}(t_{n-1},\cdot)$ by applying the regression based  DNN algorithm method from Section~\ref{sec:linearPIDE}. In this step it is important to work with  DNNs with a smooth activation function so that $\nabla {\mathcal U}_n$  is well-defined.\footnote {If $g(T,\cdot)$ is not $\mathcal{C}^1$  one has to define $\mathcal{U}_N$ as  DNN approximation to the terminal condition; see for instance~\citet{germain2020deep}.}
%

\begin{remark} \label{remark:trunc}
An additional complication may arise if the  domain $D$ is unbounded. In that case we compute the approximate solution on some bounded set $A_n \subset D$ and we have to make sure that the set $A_n$ is big enough so that the probability of generating  paths with $X_{t_{n-1}} \in A_{n-1}$ but $X_{t_{n}} \notin A_{n}$ is sufficiently small so that these paths can be ignored without affecting the training procedure.  There are various ways to achieve this, see \citet{beck2019deep} for details.
\end{remark}

\subsection{Alternative method} Next we discuss an alternative linearization procedure based on the \emph{midpoint rule}
 \begin{align} \label{eq:midpoint1}
\int_{t_{n-1}}^{t_{n}} f_{[u^*]}(s,X^x_s) \ud s \approx  f_{[u^*]}\left(\overline {t}_{n},X^x_{\overline{t}_n})\right)  \Delta t_n\,.
\end{align}
Here $ \overline t_n = (t_{n-1} + t_{n})/2$ is the midpoint of the interval $[t_{n-1}, t_{n}]$. It is well known that \eqref{eq:midpoint1} usually provides  a better approximation to an integral than the endpoint rule \eqref{eq:endpoint1}.
We therefore propose the following  approximation scheme based on the midpoint rule. In a first step we apply the approximation~\eqref{eq:basic-linearization}  over the smaller interval $[\overline t_n, t_{n}]$; this yields an approximate solution $\overline{\mathcal{U}}_{n}$ at $\overline{t}_n$. Next we define
\begin{align}\label{eq:semilinear_optim}
	\tilde{u}(t_{n-1},x) &:= \mathbb{E} \bigg[{\mathcal U}_n(X^x_{t_n }) \I_{\{\tau >t_n\}} + g(\tau, X^x_\tau)   \I_{\{\tau \le t_n\}} \\&+  (t_n \wedge \tau - t_{n-1})  f\big (\overline{t}_n \wedge \tau ,  X^x_{\overline{t}_n \wedge \tau}, \overline{\mathcal U}_{n} ( X^x_{\overline{t}_n \wedge \tau}), \nabla \overline{\mathcal U}_{ n} ( X^x_{\overline{t}_n \wedge \tau})  \big)  \bigg]\,;
\end{align}
as before, a numerical approximation $\mathcal{U}_{n-1}$ to $\tilde{u}(t_{n-1},x)$ is computed  via the deep learning algorithm from Section~\ref{sec:linearPIDE}.
Our numerical experiments in Section~\ref{sec:examples-semilinear} show that  the approximation based on the midpoint rule performs significantly better than the original deep splitting  method from \citet{beck2019deep}.

\begin{remark}[Convergence results] \citet{germain2020deep}   recently obtained  a convergence result for the  deep splitting method for PDEs. They consider only Cauchy problems (no boundary conditions) and they make strong Lipschitz assumptions on the coefficients of the equation. Under these conditions  they come up  with  an estimate for the approximation error in terms of the mesh of the partition $t_0, \dots, t_N$ and  the size of the DNN used in the individual time steps.   We are confident that similar results can be obtained for PIDEs. However, theoretical convergence  results are of limited practical use for designing an effective DNN  approximation to PIDEs,  and we leave this  technical issue for future work.
\end{remark}

\section{Examples for the semilinear case} \label{sec:examples-semilinear}

In this section we test  the algorithm in two case studies. First we consider the well-known  stochastic regulator problem. This problem leads to a semilinear PIDE with an explicit solution  that can be used  as a validity check for our methods. The second case study  is an actuarial example dealing with the optimization of an insurance portfolio  under transaction costs and risk capital constraints.

\subsection{Stochastic regulator problem}
The first example is the \emph{stochastic linear regulator}. Denote by $c =(c_{1,t}, \dots, c_{d,t})_{t\ge 0}$  an adapted  control strategy and consider the controlled $d$-dimensional process $\widetilde X^c$ with dynamics
\begin{align}\label{eq:state-process-regulator}
d\widetilde X_{i,t}=c_{i,t} \ \ud t + \sigma_i  \ud W_{i,t} + \int_{\R}z \ 	\widetilde{J_i}(\ud z,\ud t), \; 1 \le i \le d,  \qquad \widetilde X_0=x \in \R^d\, .	
\end{align}
Here $W = (W_1, \dots, W_d)$ is a $d$-dimensional standard Brownian motion,  $\theta \in \R^d$, $\rho \in \R^d $, $\sigma_1, \dots, \sigma_d$ are positive  constants, $T >0$   and $\widetilde J_1, \dots, \widetilde{J}_d$ is the compensated jump measure of $d$ independent compound   Poisson processes with Gamma($\alpha_i,\beta_i$)-distributed jumps. Denote by $\mathcal{A}$ the set of all adapted $d$-dimensional processes $c$ with $\bE \Big [ \int_0^T |c_s|^2 \ud s\Big] < \infty $  and  consider the control problem
\begin{align*}
u(t,x)=\inf_{c \in \mathcal{A}} \bE \bigg [ \sum_{i=1}^d \Big( \int_{t}^{T} \big( (\widetilde X_{i,s}^c)^2 + \theta_i c_{i,s}^2 \big) \,\ud s + \rho_i (\widetilde X_{i,T}^c )^2 \Big) \,\Big |\,  \widetilde X_t^c=x \bigg ], \quad t \in [0,T], \, x \in \R^d.
\end{align*}
 The interpretation of this problem is that the controller wants to drive  the process $\widetilde X$ to zero using the control $c$; the instantaneous control cost  (for instance the energy consumed) is  measured by $\theta  c^2_t$. At maturity $T$ the controller incurs the terminal cost $\rho \widetilde X_T^2$.

The Hamilton-Jacobi-Bellman (HJB) equation associated to this control problem  is
\begin{align*}
u_t(t,x)&+\frac{1}{2} \sum_{i=1}^{d} \sigma_i^2  u_{x_i x_i}(t,x) + \int_{\R^d} \Big [ u(t,x+z)-u(t,x)-\sum_{i=1}^d z_i u_{ x_i}(t,x) \Big ] \nu(\ud z)  \\ &+ \sum_{i=1}^d x_i^2 +\inf_c \Big \{\sum_{i=1}^d c_i u_{ x_i}(t,x) +\theta_i c_i^2  \Big \}=0,\quad (t,x) \in [0,T) \times \R^d,
\end{align*}
with terminal condition  $u(T,x)=\varphi(x):=\sum_{i=1}^d \rho_i x_i^2$.
The minimum  in the HJB equation is attained at
$c_{i}^*(t,x)=-\frac{1}{2 \theta_i} \frac{\partial u}{\partial x_i} (t,x)$, so that the value function solves the semilinear PIDE
\begin{align} \label{eq:HJB-stoch-reg}
u_t(t,x)&+\frac{1}{2} \sum_{i=1}^{d} \sigma_i^2  u_{x_i x_i}(t,x) - \sum_{i=1}^d   \int_{\R^d} z_i \ \nu(\ud z)   u_{x_i} (t,x)\\ &+\int_{\R^d}  [u(t,x+z)-u(t,x) ] \ \nu(\ud z)  + \sum_{i=1}^d x_i^2 - \sum_{i=1}^d \frac{1}{4 \theta_i} u_{ x_i}(t,x)^2=0\,.
\end{align}
It is well known that the HJB equation \eqref{eq:HJB-stoch-reg} can be solved analytically, see also \cite{oksendal2007applied}. For this we  make the Ansatz
$
u(t,x)=\sum_{i=1}^d a^i(t) x_i^2 +b(t).
$
Substitution into \eqref{eq:HJB-stoch-reg} gives an ODE system for $a(t)$ and $b(t)$ that can be solved explicitly. One obtains
\begin{align*}
			a^i(t)&=\sqrt{\theta_i} \frac{1+\kappa_i e^{2t/\sqrt{\theta}}}{1-\kappa_i e^{2t/\sqrt{\theta}}}, \qquad \kappa_i:= \frac{\rho_i-\sqrt{\theta_i}}{\rho_i+\sqrt{\theta_i}}e^{-\frac{2T}{\sqrt{\theta_i}}} \\
			b(t)&=\sum_{i=1}^d \sqrt{\theta_i} \Big(\sigma_i^2 + \int_{\R^d} z_i^2 \ \nu(\ud z)\Big)\Big((T-t)+\log \big((1-\kappa_i e^{2t})/(1-\kappa_i e^{2T})\big)\Big).
\end{align*}

\begin{figure}[h]
	\centering
	\begin{tabular}{@{}c@{\hspace{.5cm}}c@{}}
		\includegraphics[width=0.95\textwidth]{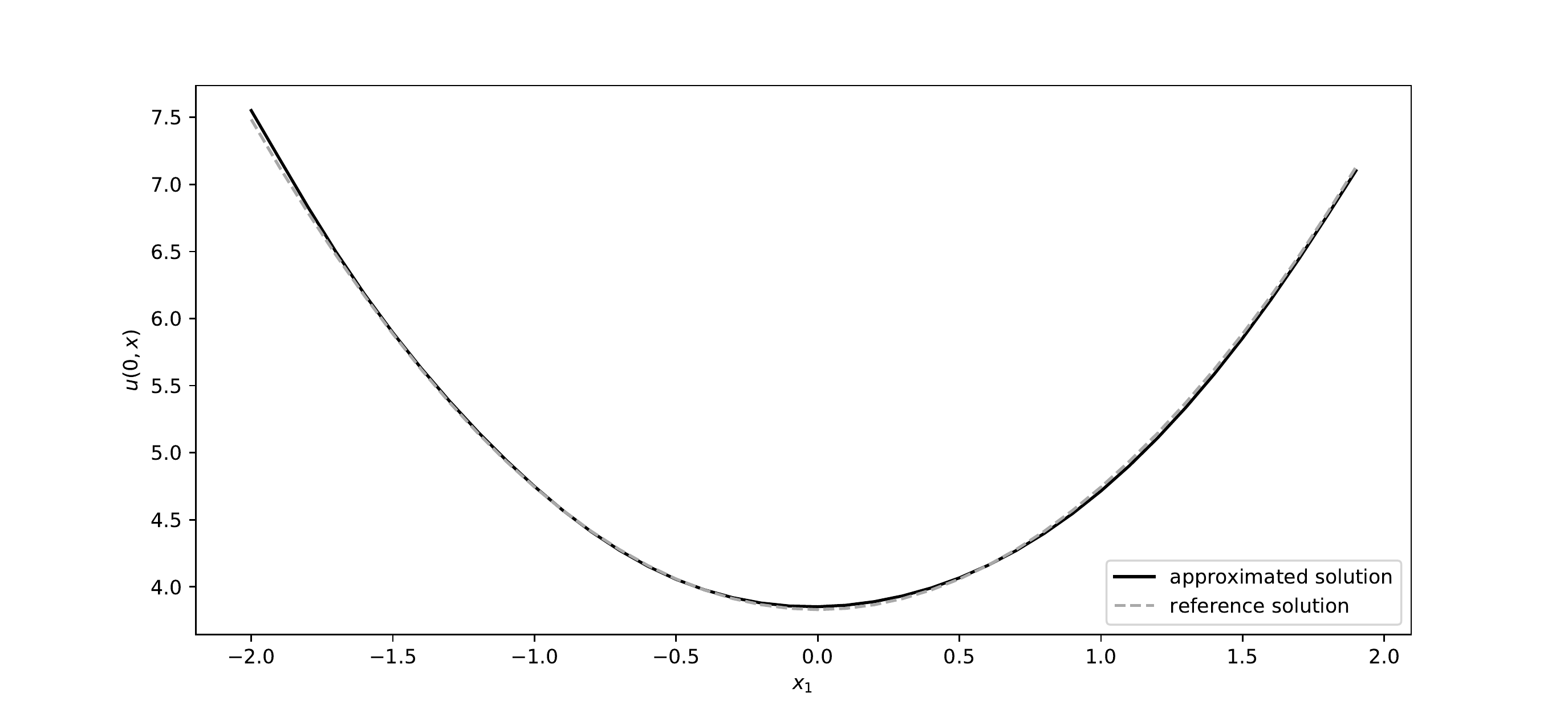}
	\end{tabular}
	\caption{ Explicit solution $u(0,x)$ (\textit{grey}) and DNN-approximation $\mathcal{U}(x)$ (\textit{black})  on  $\{x=(x_1,1,1,1): x_1 \in [-2,2]\}$.}
	\label{fig:stochreg}
\end{figure}

To test the deep splitting method we  computed $u(t,x)$ for $x \in A: =  [-2,2] ^d$ numerically. For this we partition the time horizon into $N=10$ intervals $0=t_0<t_1<\dots<t_N=T$  and simulate the auxiliary process $ X$ for $t \in [t_{n-1},t_n]$
\begin{align} \label{ex:X}
X_{i,t}=\xi_i  + \int_{t_{n-1}}^{t} \sigma_i \ \ud W_{i,s} + \int_{t_{n-1}}^{t}  \int_{\R} z \ J_i(\ud z,\ud s) - \int_{t_{n-1}}^{t}\Big( \int_{\R}z \ \nu_i(\ud z) \Big) \ \ud s, \quad 1 \le i \le d,
\end{align}
where  $W$ and $J$ are as in \eqref{eq:state-process-regulator}, and where $\xi$ is uniformly distributed on $A$.  The nonlinear term is finally given by
\begin{align*}
f(t,x,y,z)=\sum_{i=1}^d \Big( x_i^2 - \frac{1}{4\theta_i}z_i^2 \Big).
\end{align*}
We use the {midpoint rule} \eqref{eq:midpoint1} to linearize the PIDE and approximate $x \mapsto u(t_n,x)$ with a deep neural network $\mathcal{U}_{n}$ for $n=0,1,\dots,N-1$; details on numerics are given in Remark~\ref{rem:num_detail2}. Figure \ref{fig:stochreg}
shows the approximate solution of $u$ obtained by the deep splitting method and the analytic reference solution. Table \ref{tab:stochreg} contains the parameters used.

\begin{figure}[h]
	\centering
	\begin{tabular}{@{}c@{\hspace{.5cm}}c@{}}
		\includegraphics[width=0.95\textwidth]{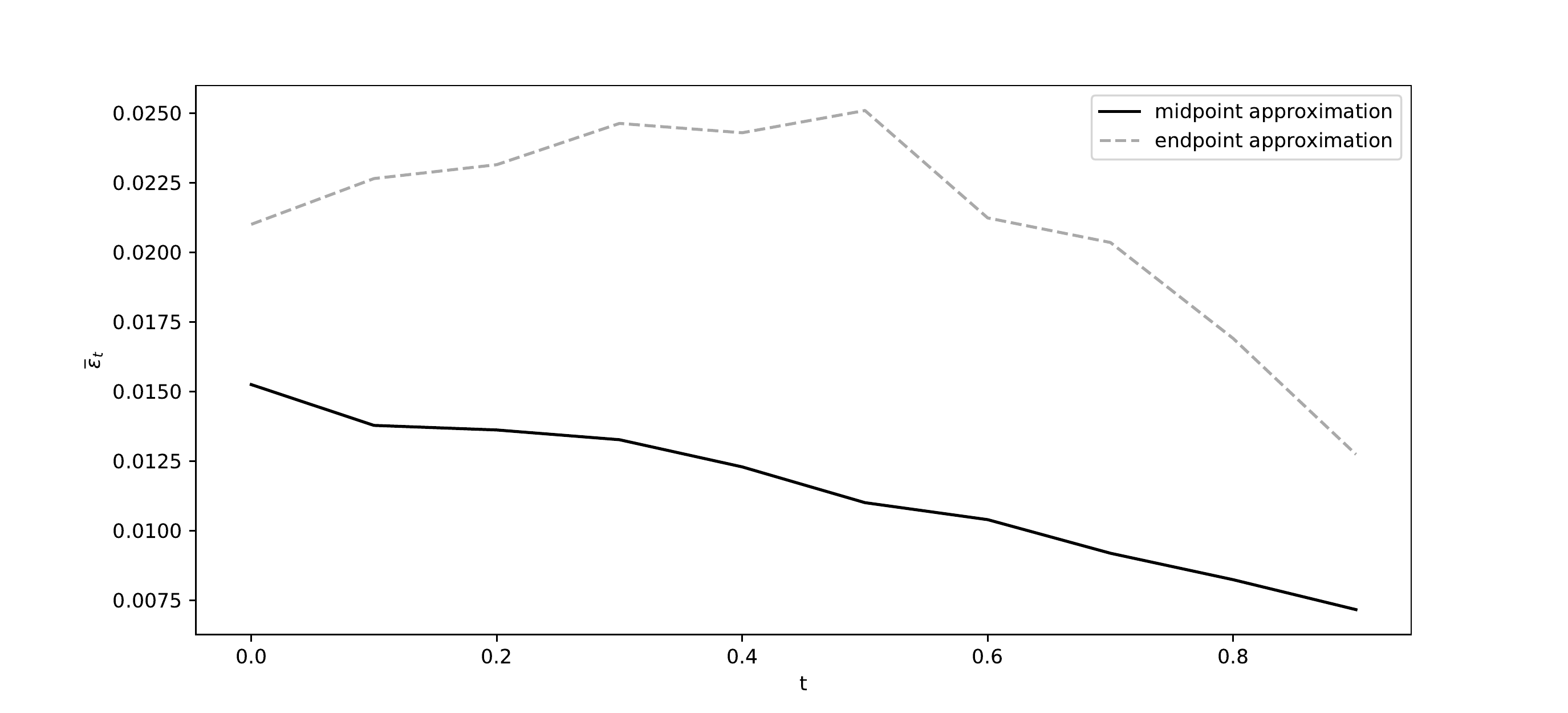}
	\end{tabular}
	\caption{Comparison of the averaged relative $L^1$-error $\bar \epsilon_t$ of the endpoint (\textit{grey}) and the midpoint (\textit{black}) approximation for $N=10$.}
	\label{fig:stoch_reg_error}
\end{figure}

As in the linear case the relative $L^1$-error for a deep neural network approximation $\mathcal{U}_{n}$ at time $t=t_n$ is defined as
\begin{align} \label{eq:rel_err_semi}
\epsilon_{t_n}:=\mathbb{E} \bigg[
\Big|\frac{\mathcal{U}_{n}(\xi)-u(t_n,\xi)}{u(t_n,\xi)}\Big|\bigg].
\end{align}
Using 10000 simulations of $\xi$ we approximate the relative error $\epsilon_t$ \eqref{eq:rel_err_semi} for $t=t_0,t_1,\dots,t_{N-1}$. In particular, for $i=1,\dots,10$ and $n=0,1,\dots,{N-1}$ we train different networks $\mathcal{U}^i_{n}$ with error $\epsilon^i_{t_n}$ to compute the average error, i.e. $\bar{\epsilon_t}=\frac{1}{10}\sum_{i =1}^{10}\epsilon^i_t$ for $t=t_n, \ n=0,1,\dots,{N-1}$.  The training for the whole network, that consists of 20 sub-networks (10 DNNs at $t_n$ and 10 DNNs at midpoints $\bar t_n$), took around 7500 seconds.
Figure \ref{fig:stoch_reg_error} compares the error $\bar \epsilon$ of the midpoint method to the averaged relative $L^1$-error computed with the {endpoint rule} \eqref{eq:endpoint1} with $N=10$. Clearly, the  method based on the {midpoint} rule performs substantially better than the standard method based on  the   endpoint rule.
Note that even though we use $N=10$ for both methods the midpoint rule requires training of 20 networks. However our experiments showed that  the midpoint rule with $N=10$ discretization steps also outperforms the endpoint rule with $N=20$ steps.

\begin{table}[h]
	\begin{center}
		\begin{tabular}{c c c c c c c c c}
			$T$ & $d$& $\sigma$& $\theta $ &  $\rho$ & $\lambda$& $\alpha$ & $\beta$ & \\
			\hline
			\rule{0pt}{0.9\normalbaselineskip}1& 4&$(0.1,\dots,0.1)$ & $(1,\dots,1)$ &  $(0.5,\dots,0.5)$ & $(10,\dots,10)$ & $(0.4,\dots,0.4)$ & $(4,\dots,4)$\\
			\hline
		\end{tabular}
	\end{center}
	\vspace{2mm}
	\caption{Parameters used in the stochastic regulator problem}  \label{tab:stochreg}
\end{table}

\begin{remark}[Details on the numeric implementation] \label{rem:num_detail2}   We used a similar network  architecture as in \cite{beck2019deep}.  We worked with  deep neural networks with 2 hidden layers consisting of 100 nodes each. All neural networks are initialized with random numbers using the Xavier initialization. We use mini-batch optimization with mini-batch size $M=5000$, batch normalization and 12000 training epochs. The loss function is minimized with the Adam optimizer and decreasing learning rate from 0.1 to 0.01, 0.001, 0.0001 and 0.00001 after 2000, 4000, 6000 and 8000 epochs. We use \texttt{sigmoid} as activation function for the hidden layers and the \texttt{identity} for the output layer. 	All computations were run on a Lenovo Thinkpad notebook with an Intel Core i5 processor (1.7 GHz) and 16 GB memory.
\end{remark}

\subsection{Optimal insurance portfolios with risk capital  constraints}

The second case study is an optimization problem in  insurance. We consider a stylized model of an insurance company  who dynamically optimizes  her holdings of  some insurance portfolio in the presence of transaction costs and risk capital constraints.   Without of  capital constraints the problem admits  an analytic solution and  is therefore a useful test case. With capital constraints the model leads to a semilinear boundary value problem for which there is no explicit solution, and we study this case numerically using the deep splitting method. In fact, the analysis of the impact of risk capital constraints on optimal insurance portfolios is also interesting in its own right.
		
Consider an insurer who  invests into an insurance portfolio with market value $S_t$ at time $t$ and in cash. We assume that the process $S$ has dynamics
$$
\ud S_t =\bar{\mu} \ud t + \sigma \ud W_t - \ud R_t,
$$
Here the Brownian motion $W$ models the compensated small claims, the compound Poisson process $R$ with intensity $\lambda$ and jump size distribution $\eta(\ud z)$ models  the large claims, and $\bar{\mu} >0$ is the  premium income of the portfolio.  We assume that the insurance portfolio can be adjusted only gradually over time. In mathematical terms this means that the insurer uses an absolutely continuous  trading strategy  with  trading rate $\theta= (\theta_t)_{0\le t \le T}$, so that for a given strategy $\theta$, the size $Q^\theta$ of the insurance portfolio  satisfies $\ud  Q_t^\theta =\theta_t \ \ud t$. We assume that
Moreover,  the insurer incurs a transaction cost that is proportional to the trading rate $\theta_t^2$.   This   models the fact that a rapid  adjustment of the insurance portfolio is expensive, as the insurer needs to carry out costly marketing activities  or to enter into expensive reinsurance contracts. Given a parameter $\kappa >0$ modelling the size of the transaction cost, the cash position $C^\theta$ has dynamics
$$				\ud C_t^\theta = -\big ( \theta_t S_t + \kappa \theta_t^2\big ) \ud t\,. $$
Here the  term  $-\theta_t S_t \ud t$  gives the  revenue from trading the insurance position and the term $- \kappa \theta_t^2 \ud t$ models the reduction of the  cash position due to transaction costs.
From a mathematical viewpoint this model is fairly similar to models used in the literature on optimal portfolio execution  such as \citet{bib:cardaliaguet-leHalle-18} or  \citet{bib:cartea-jaimungal-penalva-15}.

Denote by $E_t^\theta =S_t Q_t^\theta +C_t^\theta$ the insurer's \emph{equity}. We assume that the insurer wants to maximize  the expected value  of  her equity position $E_T^\theta$ at some horizon date time $T$, and that  she incurs a   liquidation cost of size $\gamma ({  Q_T^\theta})^2$ for some parameter $\gamma >0$.
Note that for a given strategy $\theta$ the process $(Q^\theta,E^\theta )$ has dynamics
		\begin{align} \label{statevar}
		\ud  Q_t^\theta&=\theta_t \ \ud t, \\
		\ud E_t^\theta &=  Q_t^\theta \bar{\mu} \ud t +   Q_t^\theta \sigma \ud W_t -   Q_t^\theta \ud R_t-\kappa \theta_t^2 \ \ud t. \nonumber
		\end{align}
Hence the pair $(Q^\theta, E^\theta)$ is Markov  and we may use this process as  state process.   We define the value function of the insurers optimization problem for $t \in [0,T]$ as
\begin{equation}\label{eq:opt-insurance}
		u(t,q,e):=\sup_{\theta \in \mathcal{A}} \mathbb{E}\big[E_T^\theta-\gamma (  Q_T^\theta)^2 \mid     Q_t^\theta=q , \ E_t^\theta=e \big]\,,
		\end{equation}
where $\mathcal{A}$ denotes the set of all adapted trading rates such that $E\big [\int_0^T \theta_t^2 \ud t \big ] < \infty$.

To prevent early ruin regulatory institutions usually impose risk capital constraints. In particular they may liquidate an insurance company if the  equity is too low compared to the size of the company's insurance portfolio.
Below we use the  deep splitting method to study  how  capital constraints affects the value function and the optimal strategy  of the reinsurer. Before that we consider the unconstraint problem.   In that case  there exists an explicit solution to the HJB equation (PIDE \eqref{eq:PIDE_opt-insurance} below), which  can be used  to test  the performance of the deep splitting method.

\subsubsection{The case without constraints}
	By standard arguments the HJB equation associated with the  problem \eqref{eq:opt-insurance} is
		\begin{align}\label{eq:HJB-opt-insurance}
		 u_t(t,q,e) & + q \bar{\mu} u_{e}(t,q,e) + \frac{1}{2}\sigma^2 q^2 u_{ee}(t,q,e) + \lambda  \int_{0}^{\infty}  \! [u(t,q,e+qz)-u(t,q,e)]  \eta(\ud z) \\ &+ \sup_{\theta \in \R} \{\theta u_q(t,q,e)-\kappa \theta^2 u_e(t,q,e) \}  =0,
		\end{align}
		with terminal condition $u(T,q,e)=e^2-\gamma q^2$. Maximization gives the candidate optimal trading rate
	\begin{align}\label{eq:rate-opt-insurance}
		\theta_t^*(t,q,e)=\frac{1}{2 \kappa } 		\frac{u_q(t,q,e)}{u_e(t,q,e)};
	\end{align}
		substitution into \eqref{eq:HJB-opt-insurance} yields the semilinear PIDE
\begin{align}\label{eq:PIDE_opt-insurance}
		u_t(t,q,e) &+ q \bar{\mu} u_e(t,q,e) + \frac{1}{2}\sigma^2 q^2 u_{ee}(t,q,e) + \lambda \int_{0}^{\infty} [u(t,q,e+qz)-u(t,q,e)] \ \eta(\ud z)
              \\ & + \frac{1}{4 \kappa } 		\frac{u_q(t,q,e)^2}{u_e(t,q,e)}    =0\,.
\end{align}
To solve the case without constraints   we make the Ansatz $u(t,q,e)=e + v(t,q)$, that is we assume that $u$ is linear in the equity value  $e$. This implies
 $$\lambda   \int_{0}^{\infty} [u(t,q,e+qz)-u(t,q,e)]   \eta (\ud z)=\lambda q  \int_{0}^{\infty} z   \eta (\ud z)=:q \lambda \bar \eta. $$
 If we define $\bar{\alpha}=\bar{\mu}-\lambda \bar{\eta}$, \eqref{eq:PIDE_opt-insurance} reduces to the following the following first order PDE  for $v$
\begin{equation} \label{eq:HJBv}		
-\bar \alpha q=v_t (t,q)+ \frac{v_q(t,q)^2}{4 \kappa}\,,
\end{equation}
with terminal condition $v(T,q)=-\gamma q^2$.		To find an explicit solution for $v$  we follow  \citet{bib:cardaliaguet-leHalle-18} and make the Ansatz		
$v(t,q)=h_0(t)+h_1(t)q-\frac{1}{2}h_2(t)q ^2.$
		Substitution in the HJB equation \eqref{eq:HJBv} gives  the following ODE system for $h_0, h_1, h_2$
		\begin{align*}
			h_2'=h_2^2/(2 \kappa), \qquad  h_1' =-\bar{\alpha}+h_1h_2/(2\kappa), \qquad h_0'=-h_1^2/(4\kappa)
		\end{align*}
		with terminal condition $h_0(T)=h_1(T)=0$ and $h_2(T)=2\gamma$. The ODEs can be solved explicitly: one has
		\begin{align*}
		h_2(t)&=\frac{2 \kappa}{T+\frac{\kappa}{\gamma}-t}, \\
		h_1(t)&= \frac{c_1+\frac{1}{2}\bar{\alpha}t^2-a\bar{\alpha}t}{a-t}, \text{ where }  c_1:= -\frac{1}{2}\bar{\alpha}T^2+a\bar{\alpha}T \text{ and } a:=T+\frac{\kappa}{\gamma}, \\
		h_0(t)&=\frac{1}{4\kappa (a-t)}\Big(  -\frac{2}{3}\bar \alpha^2 a^4+\frac{2}{3}\bar \alpha^2 a^3 t + 2 \bar \alpha a^2 c_1+\bar \alpha c_1 t^2 -2 \bar \alpha a c_1 t-c_1^2 + \frac{1}{12}\bar \alpha^2 t^4 -\frac{1}{3}\bar \alpha^2 a t^3 \Big).		\end{align*}

\begin{figure}[h]
			\centering
			\begin{tabular}{@{}c@{\hspace{.5cm}}c@{}}
				\includegraphics[width=0.95\textwidth]{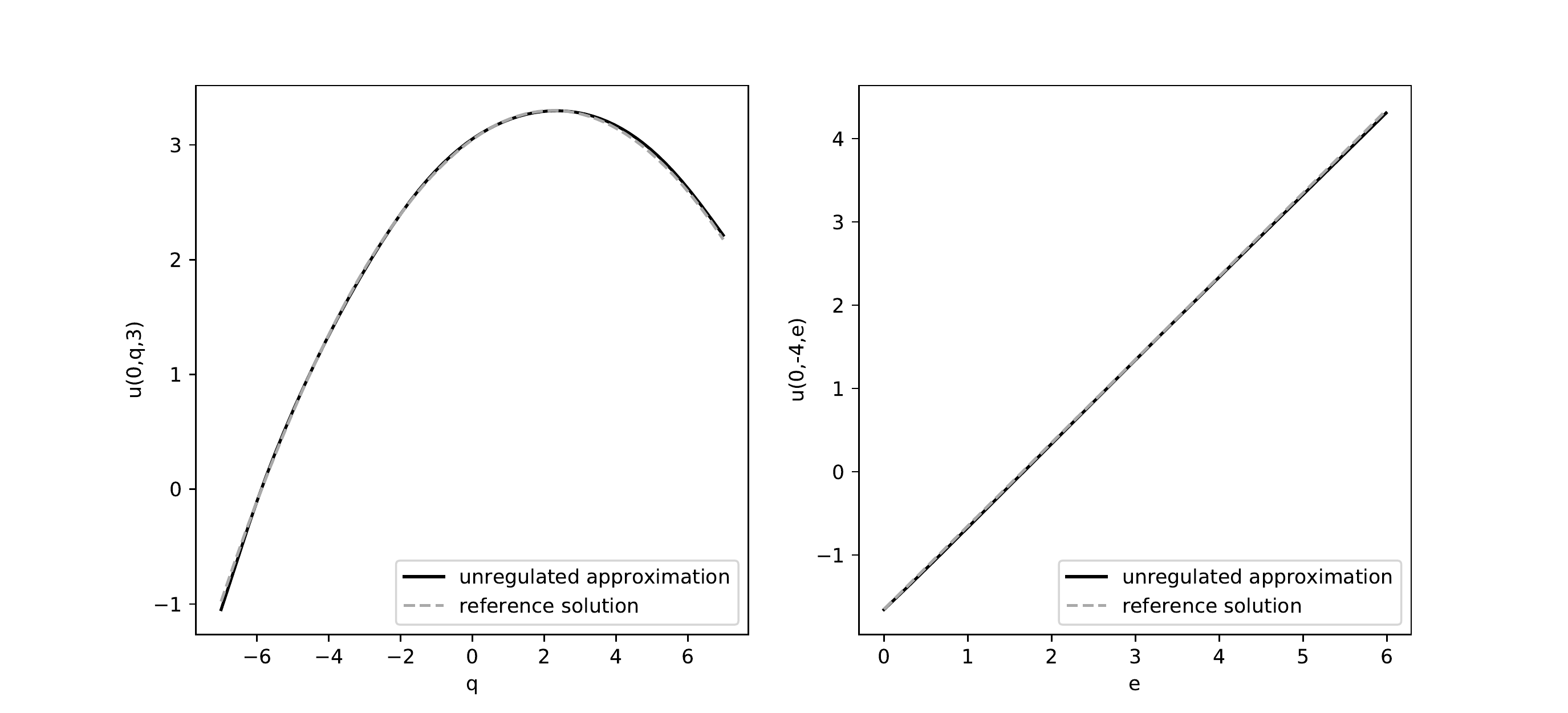}
			\end{tabular}
			\caption{ Solution of the unconstraint insurance problem for  $q \in [-7,7]$, $e=3$ (\textit{left}) and  $q=-4$, $e \in [0,6]$ (\textit{right}) with the DNN-algorithm (\textit{black}) and reference solution (\textit{grey}).}
			\label{fig:optre}
		\end{figure}

For comparison purposes we computed the solution  $u(t,x) =u(t,q,e)$ of \eqref{eq:PIDE_opt-insurance}  with the deep splitting  algorithm  on the set A:=$\{(q,e)\colon -7\leq q\leq 7, 0\leq e \leq 6\}$.    For this we partition the time horizon into $N=10$ intervals $0=t_0<t_1<\dots<t_N=T$  and simulate the auxiliary process $ X =(X^Q, X^E) $ with  $X^Q_{t}=\xi^Q$ and
	\begin{align}  \label{eq:XQ}
			 X^E_{t} &=\xi^E + \int_{t_{n-1}}^{t}	 X^Q_{s} \bar{\mu}\ \ud s + \int_{t_{n-1}}^{t}  X^Q_{s} \sigma \ \ud W_s + \int_{t_{n-1}}^{t}  \int_{\R}  X^Q_{s}  z \ J(\ud z,\ud s)\, ,
		\end{align}
where $\xi=(\xi^Q, \xi^E)$ is uniformly distributed on $[-7,7]\times [0,7]$, $W$ is a  Brownian motion and $J$ is a jump measure with compensating measure $\nu=\lambda \eta$ for $\eta$ a Gamma distribution with parameters $\alpha$ and $\beta$. Moreover, we identify the nonlinear term in the PIDE to
		\begin{align*}
		f(t,x,y,z)=\frac{1}{4 \kappa} \frac{z_1^2}{z_2}.
		\end{align*}
We used the {midpoint rule} \eqref{eq:midpoint1} to linearize the PIDE. The network architecture was as in Remark~\ref{rem:num_detail2}. The training for the whole network that consists of 20 subnetworks (a  DNN for  every grid point  $t_n$, $0\le n \le 9$ and a  DNN for every  midpoint $\bar t_n$) took around 4400 seconds.
Table \ref{tab:optre} contains the parameters used in the experiment. Figure \ref{fig:optre}
		shows the approximate solution $\mathcal{U}_0$ obtained by the deep splitting algorithm and the analytic reference solution.

		\begin{table}[h]
			\begin{center}
				\begin{tabular}{c c c c c c c  }
					$\sigma$& $\gamma $ &   $\kappa$ & $\bar \mu$ &$\lambda$& $\alpha$ & $\beta$  \\
					\hline
					\rule{0pt}{0.9\normalbaselineskip}
					$0.1$ & $0.1$ &  $0.1$ & $0.8$ & $5$ & $0.4$ & $4$\\
					\hline
				\end{tabular}
			\end{center}
			\vspace{2mm}
			\caption{Parameters for the optimal insurance problem}  \label{tab:optre}
		\end{table}
	
\begin{figure}[h]
			\centering
			\begin{tabular}{@{}c@{\hspace{.5cm}}c@{}}
				\includegraphics[width=0.95\textwidth]{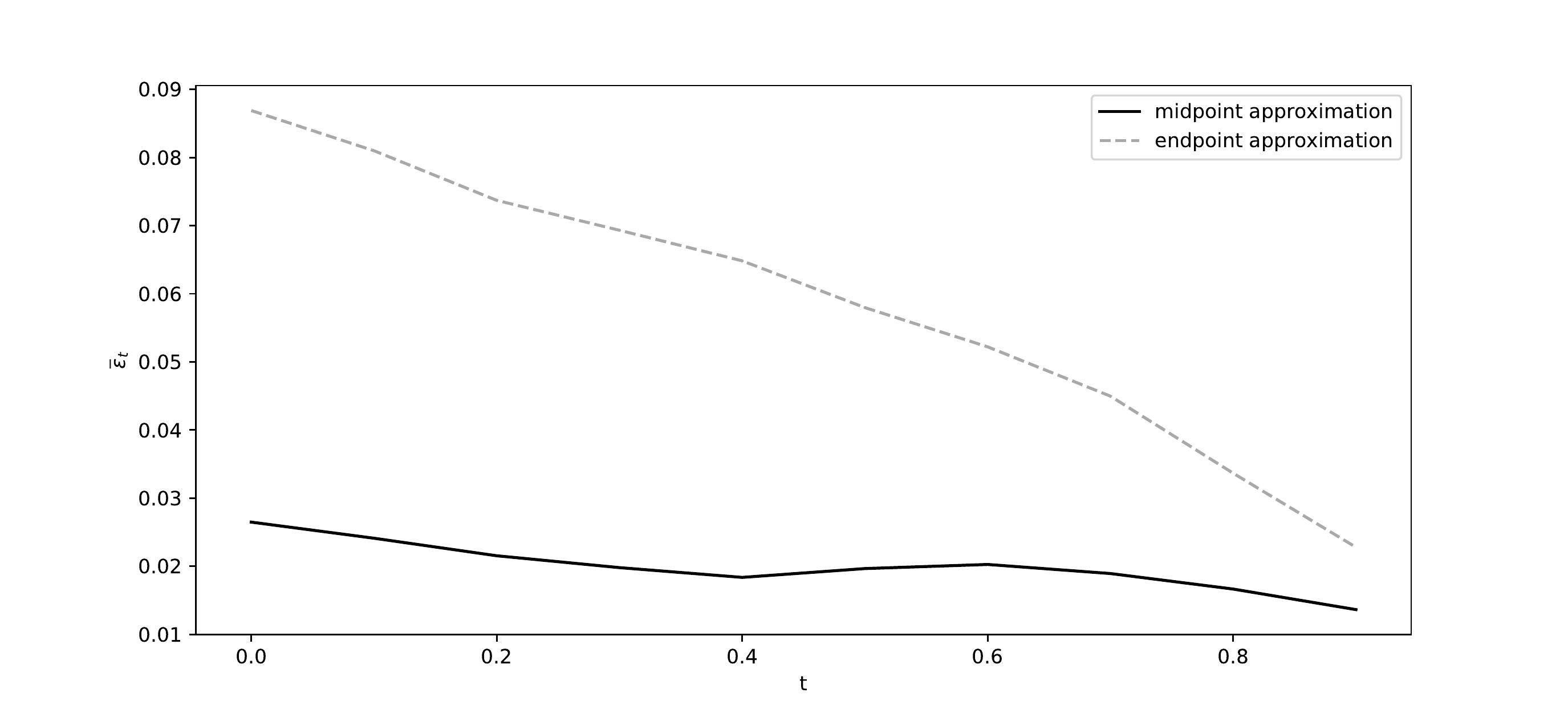}
			\end{tabular}
			\caption{Comparison of the two different training methods illustrating the averaged absolute $L^1$-error $\bar \varepsilon_t^{\text {abs}}$ using endpoint (\textit{grey}) and midpoint (\textit{black}) approximation in 10 training procedures.}
			\label{fig:error_optre}
		\end{figure}

As the solution $u$ can be zero we report  the \emph{absolute $L^1$-error} between the DNN approximation $\mathcal{U}_{n}$ and the analytic solution, that is  $		\varepsilon_{t_n}^{\text {abs}}:=\mathbb{E} \Big[\big|{\mathcal{U}_{n}(\xi)-u(t_n,\xi)}\big|\Big].$
To compute the absolute error we used  10000 simulations of $\xi=(\xi^Q,\xi^E)$ where $\xi^Q \sim \text{Unif}([-7,7])$ and  $\xi^E \sim \text{Unif}([0,6])$
		Figure \ref{fig:error_optre} illustrates the mean error of 10 different DNN approximations $\mathcal{U}^i_n, \, i=1,\dots,10$, i.e.
 using once the  {midpoint rule} and once the original method  from \cite{beck2019deep}. We see that the algorithm based on the  midpoint rule performs substantially better than the original deep splitting  method.

		\subsubsection{The optimization problem with constraints.} Next we discuss  the optimization problem under risk capital constraints.
We expect that  risk capital constraints will alter the investment decisions  of the insurer since she wants to avoid  a liquidation of the company. To set up the corresponding optimization problem we need to specify the set $D$ of acceptable positions of the insurer  and the liquidation value of the company given that its position exits  from $D$. For numerical reasons it is convenient to model  $D$ as a bounded set. Given large constants  $\bar e$ and $\bar q$  and a parameter  $\delta >0$ we define the set $D$ of acceptable positions as
\begin{align} \label{eq:acceptable-positions}
			D=\{(q,e)\in \R^2 \colon |q|<\bar q, \ \delta |q|<e<\bar e\}\,.
\end{align}
The constraint $e > \delta |q|$ implies that before liquidation  the \emph{solvency ratio}  $|Q_t| /E_t$ of the insurer is bounded by  $\delta^{-1}$, so that $\delta^{-1}$ can be viewed as a measure of the regulator's risk tolerance.   We assume that the liquidation value of the insurance company on the boundary of $D$ is given by the value $u^{\text{unreg}}$  of the company in the unconstraint problem,   reduced by a penalty $k(t,q,e)\ge 0$, and we set
$$ g(t,q,e) = u^{\text{unreg}}(t,q,e) - k(t,q,e).$$
Here $ k:[0,T]\times[-\bar q,\bar q] \times[0,\bar e] \rightarrow [0,\infty)$ is a smooth function that penalizes exit from the set of acceptable positions. We assume that
$k(t,q,\bar e)=0, 0 \leq t \leq T$ (no penalization for high equity values) and  that  $k(T,q,e)=0, \ 0 \leq e \leq \bar e$ (no penalization at maturity).
In our numerical experiments we take $\bar e = 10^6$  $\delta=1$,  and $k(t,q,e) = 0.5(T-t)(\bar e -e)/\bar e$; the other parameters are identical to the unconstraint case, see  Table \ref{tab:optre}.

The value function for the problem with risk capital constraints is
		\begin{align*}
			u^{\text{reg}}(t,q,e)=\sup_{\theta \in \mathcal{A}}
                \bE \big [g(\tau \wedge T, Q^\theta_{\tau \wedge T},E^\theta_{\tau \wedge T})|  Q_t=q, \ E_t=e \big], \quad t\in[0,T], \,  (q,e) \in D \,,
        \end{align*}
        where $\tau = \inf\{s \ge t\colon (Q^\theta,E^\theta) \notin D\}$.
The PIDE associated with this control problem and the form of the optimal strategy are  the same as in the unregulated case (see \eqref{eq:HJB-opt-insurance} and \eqref{eq:rate-opt-insurance}), but now we have in addition the  boundary condition
$$ u(t,q,e)=u^{\text{unreg}}(t,q,e)-k(t,q,e), \quad t\in[0,T), \,  (q,e) \notin D \,.$$

\begin{figure}[h]
			\centering
			\begin{tabular}{@{}c@{\hspace{.5cm}}c@{}}
				\includegraphics[width=0.95\textwidth]{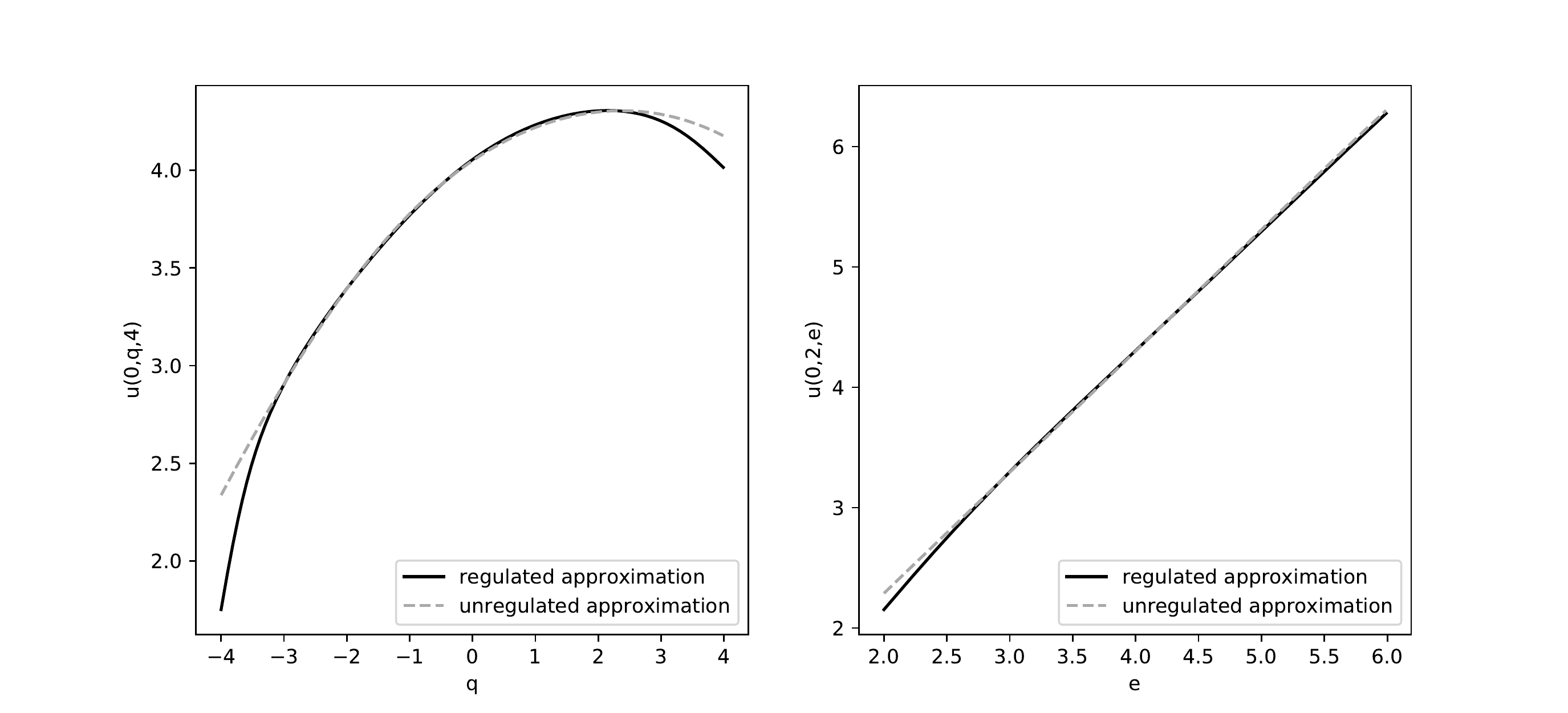}
			\end{tabular}
			\caption{ DNN-approximations of $u^{reg}(0,q,e)$ (\textit{black}) and $u^{unreg}(0,q,e)$ (\textit{grey}) on the section $\{(q,4): q \in [-4,4]\}$ (\textit{left}) and on  $\{(2,e): e \in [2,6]\}$  (\textit{right}).}
			\label{fig:optre_reg}
		\end{figure}

		Due to the boundary condition  the Ansatz 	$u^{\text{reg}}(t,q,e) = e + v(t,q)$ does not hold  and  we had to  compute $u^{\text{reg}}$  via the deep splitting algorithm.  In view of its superior performance in the unconstraint case we worked with  the \textit{midpoint} procedure, and we used the same network architecture as in the unconstraint case. We worked on the set $ A= \{(q,e) \colon -e \le q \le e, 0 \le e \leq 6 \}  \subset D$. Figure~\ref{fig:optre_reg} illustrates the DNN-approximation for $u^{\text{reg}}(0,q,e)$ and $u^{\text{unreg}}(0,q,e)$ on  the sections $\{(q,4): q \in [-4,4]\}$ and $\{(2,e): e \in [2,6]\}$. The right plot  shows that  with risk capital constraints  the value function is concave in  $e$ and for $e\approx 2$ significantly  lower than $u^{\text{unreg}}$  (for $\delta =1$ the point $(q,e)= (2,2)$ belongs to the lower bound of $D$).

\begin{figure}[h]
			\centering
			\begin{tabular}{@{}c@{\hspace{.5cm}}c@{}}
				\includegraphics[width=0.95\textwidth]{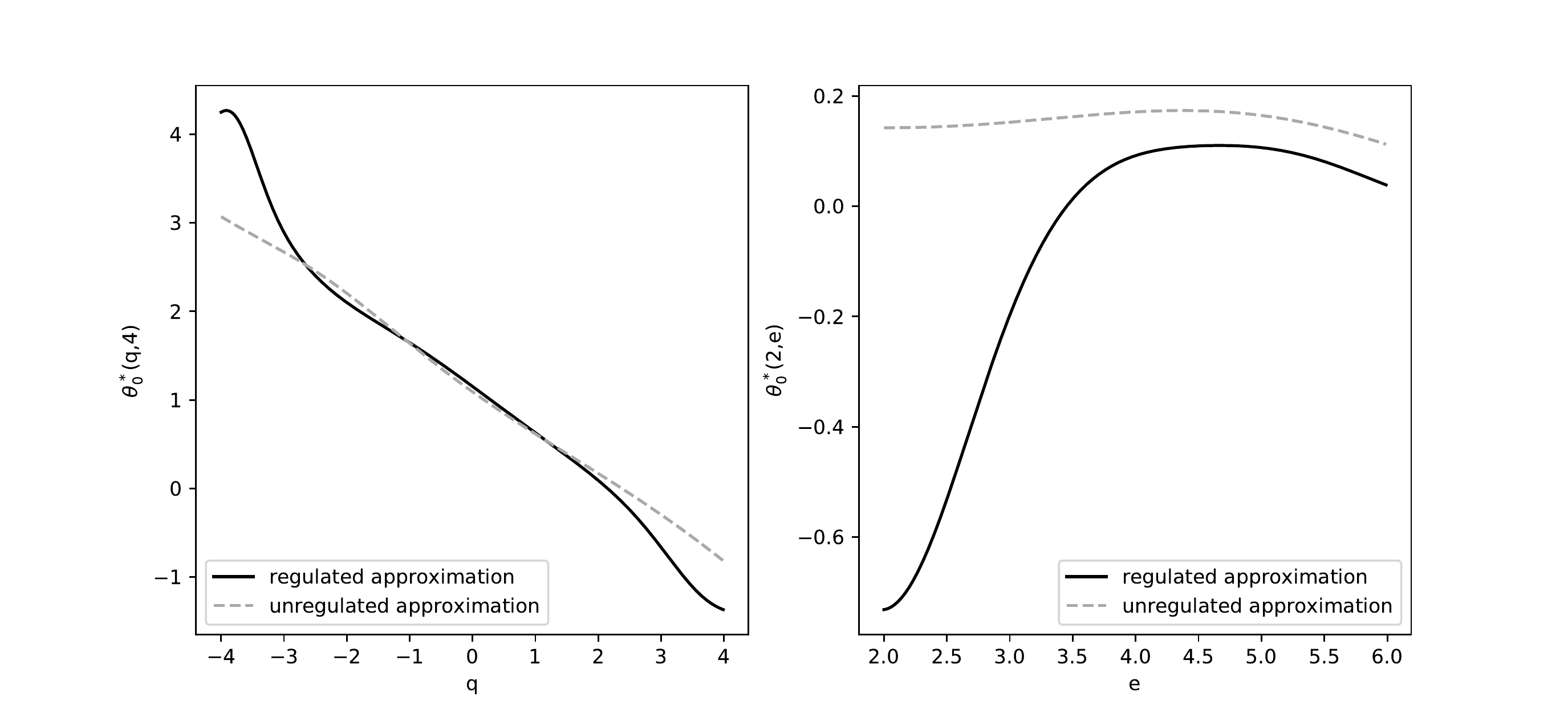}
			\end{tabular}
			\caption{ DNN-approximations of the optimal strategy $\theta_0^*$ for the regulated case (\textit{black}) and the unregulated case  (\textit{grey}) on $\{(q,4): q \in [-4,4]\}$ (\textit{left}) and on $\{(2,e): e \in [2,6]\}$  (\textit{right}).}
			\label{fig:optre_reg_strat}
		\end{figure}

The optimal trading rate is plotted in Figure~\ref{fig:optre_reg_strat}. The plots show that the optimal strategy differs from the optimal strategy in the unregulated case as the insurer wants to reduce the size $|\theta_t^*|$ of his risky position in order to avoid liquidation. For instance, in the right plot the optimal trading rate for $q=2$ is  negative for  $e$ close to 2.  This is quite intuitive:  the insurer wants to reduce her  insurance portfolio as the equity value approaches the boundary of  $D$ in order to avoid a costly liquidation.
		
\section{Conclusion}
In this paper we  studied two deep neural network algorithms for solving linear and semilinear parabolic PIDEs from insurance  mathematics. To assess the  performance of our methodology  we  carried out extensive tests for  several multi-dimensional  PIDEs arising in typical actuarial pricing and control  problems. In all test cases  the performance of the proposed approximation algorithms was quite satisfying in terms of accuracy and speed. This suggests  that  deep neural network algorithms might become a useful enhancement of the toolbox for solving many valuation and control problems in insurance that can be phrased in terms of a PIDE.

	\bibliographystyle{abbrvnat} 
	\bibliography{pide_bib}

\end{document}